\documentclass[11pt]{amsart}
\usepackage{amsmath}
\ifx\pdfoutput\undefined 
\usepackage{graphicx}
\else
\usepackage[pdftex]{graphicx}
\usepackage{epstopdf}
\fi

\newtheorem{theorem}{Theorem}[section]

\newtheorem{lemma}[theorem]{Lemma}

\theoremstyle{definition}

\theoremstyle{remark}

\numberwithin{equation}{section}

\newcommand{\inner}[2]{\langle#1,#2\rangle}
\newcommand{\PB}[2]{\{#1,#2\}}

\newcommand{\h}{\hbar}

\newcommand{\bbZ}{{\mathbb Z}}
\newcommand{\bbR}{{\mathbb R}}
\newcommand{\bbQ}{{\mathbb Q}}

\newcommand{\calO}{{\mathcal O}}
\newcommand{\calV}{{\mathcal V}}

\newcommand{\calW}{{\mathcal W}}

\newcommand{\calU}{{\mathcal U}}

\DeclareMathOperator{\tr}{Tr}
\DeclareMathOperator{\ad}{ad}

\newcommand{\norm}[1]{\left\Vert #1 \right\Vert}

\begin{document}
\openup3pt
\title{Some inverse spectral
results for semi-classical Schr\"odinger operators}
\author{V. Guillemin}\address{Department of Mathematics\\
Massachusetts Institute of Technology\\
Cambridge, MA 02139}
\email{vwg@math.mit.edu}
\thanks{V.G. supported in part by NSF grant DMS-0408993.}
\author{A. Uribe}
\address{Mathematics Department and Michigan Center for Theoretical Physics\\ University of Michigan\\Ann Arbor, MI 48109}
\email{uribe@umich.edu}
\thanks{A.U. supported in part by NSF grant DMS-0401064.}

%\footnote{Version of \today}

\maketitle
%\tableofcontents

\begin{abstract}
We consider a semi-classical Schr\"odinger operator,
$-\h^2\Delta + V(x)$.
Assuming that the potential admits a
unique global minimum and that
the eigenvalues of
the Hessian are linearly independent over $\bbQ$, we show that
the low-lying eigenvalues of the operator 
determine the Taylor series of 
the potential at the minimum.
\end{abstract}

\section{Introduction}

In this note we will report on some inverse spectral results for
the semi-classical Schr\"odinger operator,
\begin{equation}\label{1.1}
P = P(\h) = -\h^2\Delta + V(x).
\end{equation}
The potential, $V$, in (\ref{1.1}) will be assumed to be
in $C^\infty(\bbR^n)$ and have a unique non-degenerate global
minimum, $V(0)=0$, at $x=0$.  We will also assume that, 
for $\epsilon >0$ sufficiently small, $V^{-1}([0,\epsilon])$
is compact.  Then, for $\h$ sufficiently small,
the spectrum of $P$ in a small interval, $[0,\delta]$, consists
of a finite number of eigenvalues,  $E(\h)$.  In fact, by
Weyl's law
\begin{equation}\label{1.2}
\sharp\{\, E(\h)\;;\;0\leq E(\h)\leq\delta\,\} =
(2\pi\h)^{-n}\,\Bigl(
\mbox{Vol}\,\{\, 0\leq \norm{\xi}^2+V(x)\leq\delta\} + o(1)
\Bigr)
\end{equation}
(see for instance \cite{2}, chapter 9).

We will be concerned below with the question:  To what
extent do the eigenvalues, (\ref{1.2}), determine $V$?
In particular we will prove:

\begin{theorem}\label{One}
Assume $V$ is symmetric with respect to reflections about
the coordinate axes, i.e.\ for any choice of signs
\[
V(x_1,\ldots ,x_n) = V(\pm x_1,\ldots ,\pm x_n).
\]
In addition, assume that 
\[
V(x) = \sum_{i=1}^n u_i\,x^2_i + O(|x|^4),
\]
where $u_1,\ldots , u_n$ are linearly independent
over the rational numbers.  Then the family of
eigenvalues (\ref{1.2}) (with $\h$ in some neighborhood $(0,\h_0)$)
determines the Taylor series of $V$
at $x=0$.  In particular, if $V$ is real analytic,
these eigenvalues determine $V$.
\end{theorem}

Colin de Verdi\`ere and Zelditch have proved somewhat similar results
for the Dirichlet problem on convex regions in $\bbR^2$.
Namely, suppose $\Omega$ is a strictly convex region in the plane
which is real-analytic and invariant with respect to reflections
in the $x$ and $y$ axes.  Zelditch proves that for such
a region the Dirichlet spectrum,
\[
\Delta u = E_i u,\quad
0\not= u\in L^2(\Omega),\quad
u|_{\partial \Omega} = 0
\]
determines $\Omega$.  The idea of his proof is
to consider the billiard map, $B$, on the co-ball
bundle of $\partial \Omega$.  By a theorem of Anderson and
Melrose the singularities of the wave trace
\begin{equation}\label{1.3}
\sum_i \cos\sqrt{E_i} t
\end{equation}
occur at points $t=T_\gamma$ where $T_\gamma$ is the length
of a periodic billiard trajectory, $\gamma$.  In particular,
suppose that $\partial\Omega$ intersects the $y$ axis at points,
$(0,\pm a)$.  Let $p$ be the fixed point of
$B^2$ associated with the billiard trajectory, $\gamma$, which
goes from $(0,a)$ to $(0,-a)$ then reflects and goes back again.
Zelditch proves that the singluarities of (\ref{1.3}) at
$T=2a,\;4a,$ etc. determine the Birkhoff canonical form
of $B^2$ at $p$.  Suppose now that near $x=0$ the boundary
component of $\Omega$ sitting above the $x$ axis is the graph
of a function, $f$.  Because of the axial symmetries,
$f(x) = f(-x)$, and the boundary component of $\Omega$ sitting
below the $x$ axis is the graph of $-f$.  From these two
properties of $f$, Colin de Verdi\`ere deduces that the 
Birkhoff canonical form of $B^2$ at $p$ determines the Taylor
series of $f$ at $x=0$ and hence, if $\Omega$ is real analytic,
determines $\Omega$.  (Colin de Verdi\`ere
originally used this fact to
prove a somewhat weaker version of Zelditch's result, namely,
that regions of the type above are spectrally rigid.
See \cite{1} and \cite{9} for details.)

\bigskip
Our proof of Theorem \ref{One}
will involve roughly the  same kind of argument.
Namely, we willl prove:

\begin{theorem}\label{Two}  Modulo the assumption on $V$ above, the
Birkhoff canonical form of the classical Hamiltonian
\begin{equation}\label{1.4}
H(x,\xi) = \sum_{i=1}^n \xi_i^2 + V(x)
\end{equation}
at $x=0$, $\xi=0$ determines the Taylor series of $V$ at $x=0$.
\end{theorem}

We wil then make use of some recent results of
Iantchenko-Sj\"ostrand-Zworski to show that if $\psi\in C^\infty(\bbR)$
is a bump function with $\psi(s)=1$ for $s<1/2$ and $\psi(s)=0$
for $s>1$, then for $\epsilon >0$ sufficiently small, the trace
of the operator
\begin{equation}\label{1.5}
\psi\Bigl(\epsilon^{-1} P(\h)\Bigr)\,
\exp i t \h^{-1} P(\h)
\end{equation}
(which only involves the eigenvalues of $P(\h)$ in a 
neighborhood of zero)
determines the Birkhoff canonical form of $H$ at $(x,\xi)=(0,0)$.

We will review these results in the next section and prove Theorem
\ref{Two}
in \S 3.  In section 4 we indicate how to extend Theorem \ref{One}
to the
case when the configuration space is a Riemannian manifold.

\section{Microlocal Birkhoff canonical forms}

If $V$ is a potential with the properties described in Theorem
\ref{One}, then microlocally in a neighborhood of $(x,\xi)=(0,0)$
the Schr\"odinger operator (\ref{1.1}) can be conjugated by a 
unitary FIO to a rather simple ``quantum Birkhoff normal form".
More explicitly, there exist neighborhoods, $\calO_1$
and $\calO_2$, of $(x,\xi)=(0,0)$ in $\bbR^{2n}$, a canonical
transformation
\[
\kappa : \calO_1\to\calO_2,\quad \kappa(0,0) = (0,0),
\]
and a quantization of $\kappa$ by a unitary Fourier integral
operator, $\calU$, such that microlocallly on $\calO_2$
\begin{equation}\label{2.1}
\calU P \calU^{-1} = p(\h^2D_1^2+x_1^2,\ldots,\h^2D_n^2 + x_n^2,\h)
+Q_\infty(\h) + R
\end{equation}
where:
\begin{enumerate}
\item The smooth function $p$ is an $\h$-admissible symbol
admitting an asymptotic expansion
\begin{equation}
p(s_1,\ldots, s_n,h) \sim
\sum_{j=0}^\infty \h^j\, p_j(s_1,\ldots s_n)
\end{equation}
where the $p_j$ are smooth functions of $n$
variables, and the operator 
$p(\h^2D_1^2+x_1^2,\ldots\h^2D_n^2 + x_n^2,\h)$
%$p_j(\h^2D_1^2+x_1^2,\ldots , \h^2D_n^2 + x_n^2)$ denotes 
is the Weyl quantization
of $p(\xi_1^2+x_1^2,\ldots , \xi_n^2+x_n^2)$.
\item $p_0$ is of the form:
\begin{equation}\label{2.2}
p_0(w) = \sum_k u_k w_k^2 + \cdots
\end{equation}
the dots indicating quartic and higher order terms.
\item $Q_\infty$ is of order $\infty$ in $\h$,
\item The symbol of $R$ vanishes to infinite order at
$(x,\xi)=(0,0)$.
\end{enumerate}
We will give a breif sketch of how to prove this in the
next section.

\medskip
Let $T>0$ be the smallest period of the classical flow in
the region of phase space
$\{ H(x,\xi) = \norm{\xi}^2 + V(x) < \epsilon\,\}$.
Consider now the (smoothing) operator, (\ref{1.5}), with
$t\in (0,T)$.  We claim that its trace,  
\[
\tr(t,\h) = \sum_{j} e^{-it\h^{-1}E_j(\h)}\,\psi(\epsilon^{-1} E_j(\h))
\]
has an asymptotic expansion as $\h\to 0$,
\begin{equation}\label{2.3}
\tr(t,\h) \sim a_0(t)+\h a_1(t)+\cdots
\end{equation}
Indeed for each $t$
the operator $\psi(P(\h)) e^{-it\h^{-1}P(\h)}$ is an $\h$-Fourier
integral operator whose canonical relation is contained in
\[
\{\, (x,\xi,x',\xi')\;;\;(x',\xi') = f_t(x,\xi),\ H(x,\xi)\leq\epsilon\,\}
\]
where $f_t$ is the classical flow.  
The lemma of stationary phase yields the desired
expansion of the trace.

\medskip
The expansion, (\ref{2.3}), is unchanged if we replace $P$ by the
operator (\ref{2.1}).  
The trace of the resulting
operator, as an asymptotic series in $\h$, is equal
to the asymptotic expansion in powers of $\h$ of the expression
\begin{equation}\label{2.3a}
\sum_{k=(k_1,\ldots k_n)} 
\psi\Bigl(\epsilon^{-1}p((2k+1)\h\,,\,\h)\Bigr)\,
\exp\Bigl( i\h^{-1}t\,
p((2k+1)\h\,,\,\h)\Bigr).
\end{equation}
where we have set
$p((2k+1)\h\,,\,\h) = p((2k_1+1)\h,\ldots , (2k_n+1)\h\,,\,\h)$.
Since $\psi$ is identically equal to one in a neighborhood of
zero, as a formal power series in $\h$ the expansion of 
$\tr (t,\h)$ is equal to
\begin{equation}\label{2.3b}
\sum_{k=(k_1,\ldots k_n)}
\exp\Bigl( i\h^{-1}t\,\sum_{j=0}^\infty
\h^j\,p_j((2k+1)\h)\Bigr).
\end{equation}

Using ``Kronecker type" theorems of Stark \cite{3} and Fried \cite{4},
Iantchenko-Sj\"ostrand-Zworksi, \cite{7}, show that one can extract from
this expression the Taylor series of the functions,
$p_j(x_1^2+\xi_1^2,\ldots ,x_n^2+\xi_n^2)$, at $(x,\xi)=(0,0)$ and,
in particular, of the function, $p_0$, which is the ``classical" 
Birkhoff normal form of the classical Hamiltonian (\ref{1.4}).
To be more specific, they show that the coefficients of the
expansion of the trace are of the form
\begin{equation}\label{2.4}
a_l(t) =
q_l(it^{-1}\partial_\mu)\Pi_{j=1}^n\frac{1}{\sinh(t\mu_j/2)}|_{\mu = u}.
\end{equation} 
for some polynomials $q_l$.  In the setting of \cite{7} $t$ is an
integer, and the authors show that from the asymptotics of 
(\ref{2.4}) as $t\to\infty$ one can recover the coefficients of 
the polynomial, $q_l$.  (This uses heavily the ``Kronecker type" theorems
cited above.)  This is the key step in the reconstruction
process.  (For more details see \S 3 of \cite{7}.)

In the present case we only know the trace for 
$t\in (0,T)$ but, since (\ref{2.4}) is an analytic function
of $t$, we know the coefficients of the expansion for all $t$.
Therefore we can conclude:

\begin{theorem}(\cite{7})
The eigenvalues, (\ref{1.2}), determine the microlocal Birkhoff normal
form of $P$ (and, in particular, the classical Birkhoff normal
form of the Hamiltonian (\ref{1.4})).
\end{theorem}

\medskip
This result is a semi-classical version of earlier results of this
type by Guillemin \cite{5} and Zelditch \cite{8}.
(In fact, by the standard trick of reformulating high-energy
asymptotics as small $\h$ asymptotics, this result gives an
allternative, somewhat simpler, proof of these earlier results.)

\section{The proof of Theorem \ref{Two}}

\subsection{Classical Birkhoff canonical forms}
We begin by recalling the construction of the classical
Birkhoff canonical form.
Conjugating the Hamiltonian (\ref{1.4}) by the linear symplectomorphism
\[
x_i\mapsto u_i^{1/2}x_i,\quad
\xi_i\mapsto u_i^{-1/2}\xi_i,\quad i=1,\ldots n
\]
one can assume without loss of generality that
\[
H = \sum_i u_i(x_i^2+\xi_i^2) + V(x_1^2,\ldots,x_n^2)
\]
where $V(s_1,\ldots ,s_n) = O(s^2)$.  

\medskip
We will prove inductively that for $N=1,2,\ldots $
there exists a neighborhood, $\calO$, of $x=0=\xi$, and a 
canonical transformation, $\kappa: \calO\hookrightarrow \bbR^{2n}$,
$\kappa (0,0) = (0,0)$, such that
\begin{equation}\label{3.1}
\kappa^*H = \sum_{i=1}^{N} H_i + R_{N+1} + R_{N+1}',
\end{equation}
where:
\begin{enumerate}
\item[(a)] The $H_i$'s are homogeneous polynomials
of degree $2i$ of the form
\begin{equation}\label{3.1}
H_i = H_i(x_1^2+\xi_1^2,\ldots ,x_n^2+\xi_n^2),
\end{equation}
with 
\begin{equation}\label{3.3}
H_1 = \sum_{i=1}^n u_i\,(x_i^2+\xi_i^2).
\end{equation}
\item[(b)] $R_N$ is homogeneous of degree $2N$ and of the form:
%\[
%\tilde R_N = R_N + R_N'
%\]
%where $R_N'$ vanishes to order $2N+2$ at $(x,\xi)=(0,0)$ and
%$R_N$ is homogeneous of degree $2N$ in $(x,\xi)$ and of the form
\[
R_N = V_N + R_N^\sharp,
\]
where $V_N$ consists of the terms homogeneous of degree
$2N$ in the Taylor series of $V(x_1^2,\ldots x_n^2)$ at $x=0$,
and $R_N^\sharp$ is an artifact of the previous inductive steps.
\item[(c)] $R_N'$ vanishes to order $2N+2$ at the origin and 
is of the form
\[
R_N' = V-\sum_{k=2}^N V_k + S_N
\]
where $S_N$ is another artifact of the inductive process.
In addition, $R_N'$ is even.
\end{enumerate}
We will also show that this induction argument is such that one 
can read off from the $H_i$'s the first $N$ terms in the Taylor 
expansion of $V(s_1,\ldots s_n)$ at $s=0$.

\medskip
For $N=1$ in (\ref{3.1})
these assertions are true with $\kappa$ the identity,
$R_2=V_2$ consists of the quartic terms in
the Taylor series of $V$, and $R_2' = V-V_2$
(in particular $R_2^\sharp = 0 = S_2$).

\medskip
Let us suppose that these assertions are true for $N-1$ ($N\geq 2$),
so that
\[
\kappa^*H = \sum_{i=1}^{N-1} H_i + \tilde R_N
\]
with $\tilde R_N = R_N + R_N'$, as above.
%Let $R = R_N + R'$, where $R'$ vanishes to order $2N+2$
%at $x=0=\xi$, and $R_N$ is a homogeneous polynomial of degree $2N$.
%We can write this polynomial as a sum:  $V_N + R_N^\sharp$,
%where $V_N$ consists of the homogeneous terms of degree $2N$ in the
%Taylor series of $V$ at $x=0$ and $R_N^\sharp$ is an artifact
%of the previous $N-1$ steps in our induction.  
We look for a homogeneous polynomial of degree $2N$, $G_N$, such that
\begin{equation}\label{3.4}
\PB{H_1}{G_N} = R_N - H_N(x_1^2+\xi_1^2,\ldots, x_n^2+\xi_n^2)
\end{equation}
where $H_N(s_1,\ldots, s_n)$ is a homogeneous polynomial of degree
$N$ in $s$.  Introducing complex coordinates, 
$z_i = x_i + \sqrt{-1}\xi_i$, the Hamiltonian vector field
\[
\calV = \sum_i \frac{\partial H_1}{\partial \xi_i}
\frac{\partial\ }{\partial x_i} - 
\frac{\partial H_1}{\partial x_i}
\frac{\partial\ }{\partial \xi_i}
\]
becomes the vector field
\begin{equation}\label{3.5}
\sqrt{-1}\sum_i u_i\Bigl(
z_i\frac{\partial\ }{\partial z_i} -
\bar z_i \frac{\partial\ }{\partial \bar z_i} 
\Bigr)
\end{equation}
in these coordinates, and
\begin{equation}\label{3.6}
x_i^2+\xi_i^2 = z_i\bar z_i = |z_i|^2.
\end{equation}
Thus for $|\alpha|+|\beta| = 2N$,
\begin{equation}\label{3.7}
L_\calV\,( z^\alpha \bar z^\beta) = \sqrt{-1}
\Bigl(\sum_i u_i(\alpha_i-\beta_i) \Bigr)
z^\alpha \bar z^\beta .
\end{equation}

Suppose
$R_N =
\sum_{|\alpha|+|\beta| = 2N} c_{\alpha,\beta}z^\alpha \bar z^\beta$.
Letting 
\[
G = \frac{1}{\sqrt{-1}} \sum_{\alpha \not= \beta}
\frac{c_{\alpha,\beta}}{\inner{u}{\alpha-\beta}} z^\alpha \bar z^\beta
\]
(where we are using the hypothesis of the linear independence
of the $u_i$ over $\bbQ$), we get from (\ref{3.7})
\[
\PB{H_1}{G_N} = L_\calV\,G = R_N-H_N
\]
where
\[
H_N = \sum_{\alpha = \beta,\,|\alpha|=N}
c_{\alpha,\beta}z^\alpha \bar z^\beta.
\]
These arguments prove:
\begin{lemma}\label{3G}
There is a unique homogeneous polynomial, $G_N$, linear 
combination of monomials $z^\alpha\bar z^\beta$ with
$|\alpha|+|\beta|=2N$ and $\alpha\not= \beta$, such
that (\ref{3.4}) holds (and $H_N$ consists of the ``diagonal"
monomials of $R_N$.)
\end{lemma}

For future reference we make a little more
explicit the form of $G_N$.  Note that
\[
x_i^{2k_i} = \Bigl(\frac{z_i+\bar z_i}{2}\Bigr)^{2k_i} =
\Bigl(\frac{1}{2}\Bigr)^{2k_i}\binom{2k_i}{k_i} |z_i|^2 +
F_{k_i} + \overline{F_{k_i}},
\]
where
\[
F_{k_i} = \Bigl(\frac{1}{2}\Bigr)^{2k_i}
\sum_{0\leq r< k_i} \binom{2k_i}{r} z_i^{2k_i-r}
\bar z^{r}.
\]
This shows that
\[
x_1^{2k_1}\cdots x_n^{2k_n} = 
\binom{2k_1}{k_1}\cdots \binom{2k_n}{k_n}\
|z_1|^{2k_1}\cdots |z_n|^{2k_n} + \cdots
\]
where the  dots are a linear combination of monomials of
the form $z^\alpha \bar z^\beta$ with $\alpha\not=\beta$.
Hence by (\ref{3.6}) and (\ref{3.7}) there exists a homogeneous
polynomial,
$ G_k = G_{k_1,\ldots ,k_n}$ of degree $2N$ such that
\begin{equation}\label{3.mono}
x_1^{2k_1}\cdots x_n^{2k_n} = 
\binom{2k_1}{k_1}\cdots \binom{2k_n}{k_n}\
|z_1|^{2k_1}\cdots |z_n|^{2k_n} + \PB{H_1}{G_k}.
\end{equation}
As we'll see below, this implies that,
in solving (\ref{3.4}), we can keep track of the Taylor
coefficients of $V_N$ as terms in $H_N$ which are {\em not}
artifacts of the previous steps in our induction.

Let $\calW$ be the Hamiltonian vector field
\[
\calW = \sum_i \frac{\partial G}{\partial \xi_i}
\frac{\partial\ }{\partial x_i} -
\frac{\partial G}{\partial x_i}
\frac{\partial\ }{\partial \xi_i}.
\]
We will use the fact that, for any homogeneous function $F$, the Taylor
series of the pull-back $(\exp \calW)^*(F)$ at the origin 
is given by the expansion
\[
(\exp \calW)^*(F) \sim \sum_{k=0}^\infty (\ad_G)^k(F),
\]
where $\ad_G(F) = \PB{G}{F}$.  (Note that if $G$ is homogeneous
of degree $l$ and $F$ is homogeneous of degree $m$, then
$\PB{G}{F}$ is homogeneous of degree $m+l-2$.)
With this in mind, by (\ref{3.1}) 
\[
(\exp \calW)^*\kappa^*H = \kappa^*H + \PB{G}{\kappa^*H}+\cdots
= \kappa^*H + \PB{G}{H_1}+\cdots
\]
\[
= \kappa^*H + H_N-R_N + \cdots =
\sum_{i=1}^N H_N + \cdots
\]
where the dots represent terms that vanish to order $2N+2$ at
$(x,\xi)=(0,0)$.  
In fact, a calculation shows that for $N\geq 3$ the sum of the 
terms homogeneous of degree $2N+2$, $R_N$, equals
\[
R_{N+1} = \PB{H_2}{G}+ (R'_N)_{2N+2}
\]
where $(R'_N)_{2N+2}$ is the sum of the terms homogeneous
of degree $2N+2$ in the Taylor expansion of $R'_N$.
(For $N=2$ one has a couple of harmless
additional terms, see below.)
Thus if we replace $\kappa$ by $\kappa\exp(\calW)$,
the inductive assumptions hold
with $N-1$ replaced with $N$.  
%Moreover, $H_N$ has,
%encoded in it, the Taylor coefficients of $V_N$.  

\medskip
If we let $N$ 
tend to infinity in (\ref{3.1}) we obtain the Birkhoff canonical form
\begin{equation}\label{3.8}
\sum_{i=1}^\infty H_i(x_1^2+\xi_1^2,\ldots ,x_n^2+\xi_n^2)
\end{equation}
for the classical Hamiltonian (\ref{1.4}).
%with the
%first $N$ terms of the Taylor series of $V$
%encoded in the first $N$ terms of this sum.

\subsection{Recovering the Taylor series of $V$}
We now prove that {\em if we are given the sequence of
functions $\{H_i\}$
we can recover the Taylor series of $V$ at the origin.}

\medskip
We begin by noticing that $H_1$ consists precisely of
the quadratic terms
in the Taylor series.  Moreover, the function $G_2$ in
the first step of the inductive procedure above satisfies
\[
\PB{H_1}{G_2} = V_2-H_2
\]
(recall that $R_2=V_2$, i.\ e.\ $R_2^\sharp = 0$).
It is clear from (\ref{3.mono}) that the information in
$V_2$ is encoded in $H_2$; expilcitly, if
\[
V_2 = \sum_{k, |k|=2} c_k x^{2k},
\]
then by (\ref{3.mono})
\[
H_2 = \sum_{k, |k|=2} c_k\binom{2k}{k} |z|^{2k}
\]
where we have left $x^{2k} =x_1^{2k_1}\cdots x_n^{2k_n}$, etc.
Therefore the quartic term, $V_2$, is determined by $H_2$.
This implies that $G_2$ is also determined (see Lemma \ref{3G}).

\medskip
If we now conjugate $H=H_1+V$ by the exponential of the 
Hamilton vector field of $G_2$, we obtain:
\[
\kappa_1^*H = H_1+R_2+R_2'+
\PB{G_2}{H_1} + \PB{G_2}{R_2} + 
\PB{G_2}{\PB{G_2}{H_1}}
+O(8) =
\]
\[
= H_1+H_2+V_3+\PB{G_2}{R_2} + 
\PB{G_2}{\PB{G_2}{H_1}}
+O(8),
\]
where $O(8)$ stands for a function that vanishes to order eight
at the origin.  (We have: $R_2' = V-V_2 = V_3+O(8)$.)
Therefore
\[
R_3 = V_3+\PB{G_2}{R_2} + \PB{G_2}{\PB{G_2}{H_1}},
\]
the last two terms constituting $R_3^\sharp$, the first artifact
of the inductive process.
Notice, however, that $R_3^\sharp$ is known to us since
$G_2$ was determined in the previous step.

\medskip
The next inductive step involves the equation
\[
\PB{H_1}{G_3} = V_3+R_3^\sharp-H_3,
\]
where one should notice that $R_3^\sharp$, and therefore its
``diagonal" monomials, are known.
Since $H_3$ is also known,
arguing exactly as before (appealing to (\ref{3.mono})),
we see that $V_3$ and $G_3$ are determined by $H_1$, $H_2$ and
$H_3$.  It is clear now that one can continue indefinitely in this
fashion.

\hfill Q.\ E.\ D.

\subsection{Appendix: The existence of a quantum Birkhoff canonical
form.}
The existence of a quantum Birkhoff canonical form for the Schr\"odinger
operator, (\ref{1.1}), can be proved by essentially the same methods.
Namely, let
$p_0(s_1,\ldots ,s_n)$ be a $C^\infty$ function with the same Taylor
series as the series, $\sum_i H_i(s)$.  Then, as we have just
seen, there exists a canonical transformation, $\kappa$, conjugating
(\ref{1.4}) to the Birkhoff canonical form,
$p_0(x_1^2+\xi_1^2,\ldots ,x_n^2+\xi_n^2)$, modulo an error term
which vanishes to infinite order at $x=\xi =0$.  let $\calU$ be a
Fourier integral operator quantizing $\kappa$.  Then
\begin{equation}\label{3.9}
\calU P\calU^{-1} = p_0(\h^2D_1^2+x_1^2,\ldots ,\h^2D_n^2+x_n^2)
+\h Q + R
\end{equation}
where the symbol of $R$ vanishes to ininite order at $x=\xi =0$.
We willl prove by induction that there exists an F.\ I.\ O., $\calU_N$,
such that microlocally near $x=\xi =0$
\begin{equation}\label{3.10}
\calU_N P\calU_N^{-1} = p_0(\h^2D_1^2+x_1^2,\ldots ,\h^2D_n^2+x_n^2)
+\h^N Q + R
\end{equation}
and the symbol of $R$ vanishes to infinite order at $x=\xi =0$.
Assuming this  assertion is true for $N-1$ let us prove it for $N$.
Let the symbol of $Q$ be of the form $q=q_0(x,\xi)+\h q'(x,\xi,\h)$.
We claim
\begin{lemma}
There exist $C^\infty$ functions, $a(x,\xi)$ and 
$p_N(s_1,\ldots ,s_n)$, such that
\begin{equation}\label{3.11}
\PB{p_0}{a} = q_0 - p_N(x_1^2+\xi_1^2,\ldots ,x_n^2+\xi_n^2) + r
\end{equation}
where $r$ vanishes to infinite order at $x=\xi =0$.
\begin{proof}
By (\ref{2.2}), $p_0 = H_1 + p_0'$ where $p_0'(s) = O(s^2)$.
Writing (\ref{3.11}) in the form
\[
\PB{H_1}{a} = q_0-p_N-\PB{p_0'}{a} + r
\]
it is clear that the homogeneous terms of degree $k$ in the Taylor
series of $a$ and $p_N$ can be determined from the 
previous terms by (\ref{3.4}).
\end{proof}
\end{lemma}

Now let $A = a^W(x,\h D)$ be the Weyl quantization of $a$, and let
$V = \exp(-\sqrt{-1}\h^N A)$.  Then
\[
V(\calU_{N-1} P\calU_{N-1}^{-1})V^{-1} = 
\calU_{N-1} P \calU_{N-1}^{-1} +
\sqrt{-1}\,\h^N\,[\calU_{N-1} P\calU_{N-1}^{-1}\,,\,A]
+\h^{N+1} Q' + R'
\]
where the symbol of $R'$ vanishes to infinite order at $x=\xi =0$.
In view of (\ref{3.10}) and (\ref{3.11}) the right hand side of this
identity can be rewritten in the form
\[
\sum_{j=1}^N \h^j\,p_j(\h^2D_1^2+x_1^2,\ldots ,\h^2D_n^2+x_n^2)
+ \h^{N+1} Q'' + R''
\]
where the symbol of $R''$ vanishes to infinite order at $x=\xi =0$.
Thus with $\calU_{N-1}$ replaced by $V\calU_{N-1}$
the identity (\ref{3.10})
is valid with $N-1$ replaced by $N$, and letting $N$ tend to infinity
one obtains the Birkhoff canonical form for $P$ described in the
previous section.

\section{An extension to the Riemannian case}
Let $M$ be an $n$-dimensional
Riemannian manifold, $\Delta$ its Laplace-Beltrami
operator, and $V:M\to\bbR$ a smooth function.  
We consider the Schr\"odigner operator $P = \h^2\Delta + V$ on
$L^2(M)$.  We assume for simplicity that the  spectrum of $P$
is discrete.  In this section we describe an extension of 
Theorem 1.1 to $P$.

\medskip
We will assume that $V$ has a unique global minimum, $m\in M$,
and that it is non-degenerate.  Then the Hessian of $V$ at $m$
is a well-defined positive-definite quadratic form on $T_{m}M$.
Using the metric, we can speak of the eigenvalues, $u_1,\ldots
,u_n$, of the Hesssian.  We assume that these eigenvalues are
linearly independent over $\bbQ$.  In particular all eigenvalues
are distinct, and therefore $T_{m}M$ splits naturally as an
orthogonal direct sum of lines (spanned by eigenvectors of the
Hessian).  Denote by $G \cong \bbZ^{2^n}$ the group of linear
transformations of $T_{m}M$, generated by the reflections $v_i\mapsto
-v_i$, where $v_1,\ldots ,v_n$ is a basis of eigenvectors of the
Hessian.

\begin{theorem} Assume that 
there exists a neighborhood, $\calU$, of $m$ such that:
\begin{enumerate}
\item $G$ acts on $\calU$ by isometries.
\item $m$ is a fixed point of the action, and the
infinitesimal action on $T_{m}$ agrees with the original
linear action of $G$.
\item $V|_\calU$ is invariant under $G$.
\end{enumerate}
Assume furthermore that $V(m) = 0$.  Then, assuming that
we know the Riemannian metric on $M$, the eigenvalues
(\ref{1.2}) determine the infinite jet of $V$ at $m$.
\end{theorem}

The proof is a straightforward generalization of the one
given above in the Euclidean case.  We will say a few words
on how to show that the Birkhoff normal form of the classical
Hamiltonian,
\[
H(x,\xi) = \norm{\xi}^2+V(x)
\]
(where $(x,\xi)\in T^*M$), determines the Taylor series of
$V$ at $m$.  

Introduce geodesic normal coordinates centered at $m$,
$x = (x_1,\ldots ,x_n)$, such that the corresponding frame
at $m$ consists of eigenvectors of the Hessian.  In these
coordinates the classical Hamiltonian takes the form
\[
H(x,\xi) = \sum_{i,j=1}^n g^{ij}(x)\xi_i\xi_j + V(x)
\]
where:
\begin{enumerate}
\item $g^{ij}(x) = \delta_{ij} + O(|x|^2)$,
\item $g^{ij}(\pm x_1,\ldots ,\pm x_n) = g^{ij}(x_1,\ldots ,x_n)$,
\item $V(\pm x_1,\ldots ,\pm x_n) = V(x_1,\ldots ,x_n)$.
\end{enumerate}
(The reason for (2) and (3) is that in the chosen coordinates
$G$ acts by: $(x_1,\ldots ,x_n)\mapsto (\pm x_1,\ldots ,\pm x_n)$.)
Proceeding as in \S 3, we can rewrite the Hamiltonian in the form
\[
H(x,\xi) = \sum_i u_i(x_i^2 + \xi_i^2) + 
\sum_{i,j=1}^n h^{ij}(x_1^2,\ldots ,x_n^2)\xi_i\xi_j + V(x_1^2,\ldots ,x_n^2)
\]
where $V(s_1,\ldots, s_n) = O(s^2)$ and
$h^{ij}(s_1,\ldots ,s_n) = O(s)$. 

The construction of the Birkhoff normal form $\sum_jH_j$
proceeds as before, except that the polynomials $R_N$ are now
of the form  
\[
R_N =
\sum_{i,j=1}^n h^{ij}_{N-1}(x_1^2,\ldots ,x_n^2)\xi_i\xi_j +
V_N(x_1^2,\ldots ,x_n^2) + R_N^\sharp
\]
where $h^{ij}_N(s_1,\ldots ,s_n)$
consists of the terms homogeneous of degree
$N$ in the Taylor series of $h^{ij}(s_1,\ldots ,x_n)$.
$H_N$ consists of the ``diagonal" monomials in $R_N$, 
as before (that is, monomials of the form
$c_\alpha z^\alpha \bar z^\alpha$
in complex coordinates $z_j=x_j+\sqrt{-1}\xi_j$).

We need to show that the sequence $\{H_j\}$ determines the Taylor
series of $V$.  This is possible because we are assuming
that we know the metric, and thefore we know the
sum $\sum_{i,j=1}^n h^{ij}_{N-1}(x_1^2,\ldots ,x_n^2)\xi_i\xi_j$
appearing in $R_N$. 

For example, $H_2$ consists of the diagonal terms in
\[
R_2 = -\frac{1}{3}\sum_{i,j,k,l}R_{ikjl}(0)x_kx_l\,\xi_i\xi_j +
V_2(x_1^2,\ldots ,x_n^2).
\]
If we subtract from $H_2$ the (known) diagonal terms in
$-\frac{1}{3}\sum_{i,j,k,l}R_{ikjl}(0)x_kx_l\,\xi_i\xi_j$,
we obtain the diagonal terms in $V_2(x_1^2,\ldots ,x_n^2)$,
which, as we have seen, determine the quartic terms in the Taylor
expansion of $V$.  The higher-degree cases are no different.


\begin{thebibliography}{99}

\bibitem{1} Y. Colin de Verdi\`ere, {\it Sur  les
longeurs des trajectoires p\'eriodiqued d'un billiard},
S\'em. Sud-Rhodanien Geom. III (1983), 122-140.

\bibitem{2} M. Dimassi and J. Sj\"ostrand, Spectral
Asymptotics in the Semi-Classical Limit. Cambridge
University Press (1999).

\bibitem{3} J.J. Duistermaat and V. Guillemin,  {\it The spectrum of
positive elliptic operators and periodic bicharacteristics.}
Inv. Math.  {\bf 29} (1975), 39-79.

\bibitem{4} D. Fried, {\it Cyclic resultants of reciprocal
polynomials}, Lecture Notes in Math. 1345,
Springer Verlag (1998), 124-128.

\bibitem{5} V. Guillemin, {\it Wave trace invariants},
Duke Math. Journal {\bf 83} (1996), 287-352.

\bibitem{6} A. Iantchenko and J. Sj\"ostrand,
{\it Birkhoff normal forms for Fourier integral
operators II}, Amer. J. Math {\bf 124} (2002), 817-850.

\bibitem{7}  A. Iantchenko, J. Sj\"ostrand and
M. Zworski,
{\it Birkhoff normal forms in semi-classical
inverse problems}, Math. Res. Lett. {\bf 9} (2002), 337--362.

\bibitem{8} S. Zelditch, {\it Wave invariants at 
elliptic closed geodesics}, Geom. Funct. Anal.
{\bf 7} (1997), 145-213.

\bibitem{9} S. Zelditch, {\it Spectral determination
of analytic bi-aximetric domains},
Geom. Funct. Anal. {\bf 10} (2000), 628-677.

\end{thebibliography}
\end{document}